\documentclass[10pt]{article}
\usepackage{mathrsfs}
\usepackage{amsfonts}
\usepackage{cite} % Make references as [1-4], not [1,2,3,4]
\usepackage{url}  % Formatting web addresses
\usepackage{ifthen}  % Conditional
\usepackage{multicol}   %Columns
\usepackage[utf8]{inputenc} %unicode support
\usepackage{amsmath}
\usepackage{amsfonts}
\usepackage{amssymb}
\usepackage{float}
\usepackage{bbm}
\usepackage{latexsym}
\usepackage{graphicx}
\usepackage[noblocks]{authblk}
\numberwithin{equation}{section}
\newtheorem{theorem}{Theorem}[section]
\newtheorem{lem}{Lemma}[section]
\newtheorem{remark}{Remark}[section]
\newtheorem{pro}{Proposition}[section]
\newtheorem{definition}{Definition}[section]

\begin{document}
\title{Non-zero sum differential games of forward-backward stochastic differential delayed equations under partial information and application}

\author{Yi Zhuang\thanks{School of Mathematics, Shandong University, Jinan, PR China.}}
\maketitle

\begin{abstract}
This paper is concerned with a non-zero sum differential game problem of an anticipated forward-backward stochastic differential delayed equation under partial information. We establish a necessary maximum principle and sufficient verification theorem of the game system by virtue of the duality and convex variational method. We apply the theoretical results and stochastic filtering theory to study a linear-quadratic game system and derive the explicit form of the Nash equilibrium point and discuss the existence and uniqueness in particular cases. As an application, we consider a time-delayed pension fund manage problem with nonlinear expectation and obtain the Nash equilibrium point.
\end{abstract}


\section{Introduction}

The general nonlinear backward stochastic differential equations (BSDEs) were first developed by Pardoux and Peng \cite{PardouxPeng90}, and have been widely applied in optimal control, stochastic games, mathematical finance and related fields. If a BSDE coupled with a forward stochastic differential equation (SDE), it is called the forward-backward stochastic differential equation (FBSDE). In stochastic control area, the form of the classical Hamiltonian system is one of the FBSDEs. The classical Black-Scholes option pricing formula in the financial market can be deduced by certain FBSDE. Systems based on BSDEs or FBSDEs have been widely surveyed by many authors, see Peng \cite{Peng93,Peng97}, Karoui, Peng, and Quenez \cite{KarouiPengQuenez97}, and Yong \cite{MaYong99}, etc.

In classical case, there are many phenomena that have the nature of past-dependence, i.e. their behavior not only depends on the situation at the present time, but also on their past history. Such models were identified as  stochastic differential delayed equations (SDDEs), which  are a natural generalization of the classical SDEs and have been widely studied in engineering, life science, finance, and other fields (see for example, Mohammed \cite{Monhammed98}, Arriojas, Hu, Monhammed, and Pap \cite{AHMP07}). Recently, Chen and Wu \cite{ChenWu10} studied a stochastic control problem based on SDDE. When introducing the adjoint equation, they need some new typed of BSDEs, which had been introduced by Peng and Yang \cite{PengYang09} for the general nonlinear case and called anticipated BSDEs (ABSDEs), which also play an important role in finance and insurance (see e.g. Delong\cite{Delong12}). Moreover, a class of BSDEs with time-delayed generators (BSDDEs) has also been studied (see Wu and Wang \cite{WuWang15}, Shi and Wang \cite{ShiWang15}, Wu and Shu \cite{WuShu17}). In addition, Chen and Wu \cite{ChenWu11}, Huang, Li, and Shi \cite{HuangLiShi12} studied a linear quadratic (LQ) case based on a coupled SDDE and ABSDE called the anticipated forward-backward stochastic differential delayed equation (AFBSDDE).

Game theory has been pervading the economic theory, attracts more and more research attentions. Game theory was firstly introduced by Von Neumann and Morgenstern \cite{VonMor44}. Nash \cite{Nash51} made the fundamental contribution in Non-cooperate Games and gave the classical notion of Nash equilibrium point. Recent years, many articles on stochastic differential game problems driven by stochastic differential equations appeared. Researchers try to consider the strategy on multiple players rather than one player and try to find an equilibrium point rather than an optimal control. These problems are more complex than the classical control problems but much closer to social and behavior science. Yu \cite{Yu12lq} solved the LQ game problem on forward and backward system. \O ksendal and Sulem \cite{OksendalSulem12}, Hui and Xiao \cite{HuiXiao12} made a research on the maximum principle of forward-backward system. Chen and Yu \cite{ChenYu15} studied the maximum principle of a SDDE case, Shi and Wang \cite{ShiWang15}, Wu and Shu \cite{WuShu17} discussed a BSDDE case.

In reality, instead of complete information, there are many cases the controller can only obtains partial information, reflecting in mathematics that the control variable is adapted to a smaller filtration. Based on this phenomenon, Xiong and Zhou \cite{XiongZhou07} dealt with a Mean-Variance problem in financial market that the investor's optimal portfolio is only based on the stock and bond process he observed. This assumption of partial information is indeed natural in financing market. Recently, Wu and Wang \cite{WuWang15}, Wu and Shu \cite{WuShu17} also considered the partial information case.

From above discussion, we believe that the research on general AFBSDDEs and their wide applications in mathematical finance is important and fascinating. To our best knowledge, there are quite lacking in literature. Recently, Huang and Shi \cite{HuangShi12} discussed the optimal control problem based on AFBDDE system. Our work distinguished itself from above one in the following aspects. First, we study the stochastic differential game rather than the stochastic control system. Second, we study more practical cases that the available information to the players are partial. Third, we get a worthwhile results about the solution of the LQ case by filtering equation and solve a practical problem in financial market.

The rest of this paper is organised as follows. In section 2, we give some necessary notions and state some preliminary results. in section 3, we establish a necessary condition (maximum principle) and a sufficient condition (verification theorem) for the Nash equilibrium point. In section 4, we apply the theory discussed in Section 3 to study a linear-quadratic game problem and obtain a result of the existence and uniqueness of Nash equilibrium point in particular cases. In section 5, we study a financial problem and obtain an explicit equilibrium point.

\section{Preliminary results}

Throughout this article, we denote by $\mathbb{R}^k$ the $k$-dimensional Euclidean space, $\mathbb{R}^{k\times l}$ the collection of $k\times l$ matrices. For a given Euclidean space, we denote by $\langle\cdot,\cdot\rangle$(resp. $|\cdot|$) the scalar product(resp. norm). The superscript $\tau$ denotes the transpose of vectors or matrices.

Let $(\Omega,\mathcal{F},\{\mathcal{F}_{t}\}_{t\geq0},\mathbb{P})$ be a complete filtered probability space equipped with a $d+\bar{d}$-dimensional, $\mathcal{F}_t$-adapted standard Brownian motion $(W(\cdot),\bar{W}(\cdot))$, where $\mathcal{F}=\mathcal{F}_T$. $\mathbb{E}^{\mathcal{F}_t}[\cdot]=\mathbb{E}[\cdot|\mathcal{F}_t]$ denotes the conditional expectation under natural filtration $\mathcal{F}_t$ and $f_x(\cdot)$ denotes the partial derivative of function $f(\cdot)$ with respect to $x$. Let $T>0$ be the finite time duration and $0<\delta<T$ be the constant time delay. Moreover, we denote by $\mathbb{C}([-\delta,0];\mathbb{R}^k)$ the space of uniformly bounded continuous function on $[-\delta,0]$, by $\mathbb{L}^p_\mathcal{F}(\Omega;\mathbb{R}^k)$ the space of $\mathcal{F}$-measurable random variable $\xi$ satisfying $\mathbb{E}|\xi|^p<\infty$ for any $p\geq 1$, and by $\mathbb{L}^p_\mathcal{F}(r,s;\mathbb{R}^k)$ the space of $\mathbb{R}^k$-valued $\mathcal{F}_t$-adapted processes $\varphi(\cdot)$ satisfying $\mathbb{E}\int_r^s|\varphi(t)|^pdt<\infty$) for any $p\geq 1$.

We consider the following AFBSDDE:

\begin{equation}\label{afbsdde}
\left\{
\begin{aligned}
dx^v(t)=\ &b(t,x^v(t),x^v_\delta(t),v_{1}(t),v_{2}(t))dt+\sigma(t,x^v(t),x^v_\delta(t),v_{1}(t),v_{2}(t))dW(t)\\
&+\bar{\sigma}(t,x^v(t),x^v_\delta(t),v_{1}(t),v_{2}(t))d\bar{W}(t),\\
-dy^v(t)=\ &f(t,x^v(t),y^v(t),z^v(t),\bar{z}^v(t),y^v_{\delta^+}(t),v_1(t),v_2(t)))dt\\
&-z^v(t)dW(t)-\bar{z}^v(t)d\bar{W}(t),\quad t\in[0,T],\\
x^v(t)=\ &\xi(t),\quad t\in[-\delta,0],\\
y^v(T)=\ &G(x^v(T)),\ y^v(t)=\varphi(t),\quad t\in(T,T+\delta].\\
\end{aligned}
\right.
\end{equation}
Here $(x^v(t),y^v(t),z^v(t),\bar{z}^v(t)): \Omega\times[-\delta,T]\times[0,T+\delta]\times[0,T]\times[0,T]$, $b: \Omega\times[0,T]\times\mathbb{R}^n\times\mathbb{R}^n\times\mathbb{R}^{k_1}\times\mathbb{R}^{k_2}\rightarrow\mathbb{R}^n$,  $\sigma: \Omega\times[0,T]\times\mathbb{R}^n\times\mathbb{R}^n\times\mathbb{R}^{k_1}\times\mathbb{R}^{k_2}\rightarrow\mathbb{R}^{n\times d}$, $\bar{\sigma}: \Omega\times[0,T]\times\mathbb{R}^n\times\mathbb{R}^n\times\mathbb{R}^{k_1}\times\mathbb{R}^{k_2}\rightarrow\mathbb{R}^{n\times \bar{d}}$, $f: \Omega\times[0,T]\times\mathbb{R}^n\times\mathbb{R}^m\times\mathbb{R}^{m\times d}\times\mathbb{R}^{m\times \bar{d}}\times\mathbb{R}^m\times\mathbb{R}^{k_1}\times\mathbb{R}^{k_2}\rightarrow\mathbb{R}^m$, $G:\Omega\times\mathbb{R}^n\rightarrow\mathbb{R}^n$ are given continuous maps, $x^v_\delta(t)=x^v(t-\delta)$, $y^v_{\delta^+}(t)=\mathbb{E}^{\mathcal{F}_t}[y^v(t+\delta)]$, $\xi(\cdot)\in\mathbb{C}([-\delta,0];\mathbb{R}^n)$ is the initial path of $x^v(\cdot)$, $\varphi(\cdot)\in\mathbb{L}^2_\mathcal{F}(T,T+\delta;\mathbb{R}^m)$ is the terminal path of $y^v(\cdot)$. Here for simplicity, we omit the notation of $\omega$ in each process.

Let $U_i$ be a nonempty convex subset of $\mathbb{R}^{k_i}$, $\mathcal{G}_t\subseteq\mathcal{F}_t$ a given sub-filtration which represents the information available to the controller, and $v_i(\cdot)$ be the control process of player $i\ (i=1,2)$. We denote by $\mathcal{U}_{ad}^i$ the set of $U_i$-valued $\mathcal{G}_t$-adapted control processes $v_i(\cdot)\in\mathbb{L}_{\mathcal{G}}^2(0,T;\mathbb{R}^{k_i})$ and it is called the admissible control set for player $i\ (i=1,2)$. $\mathcal{U}_{ad}=\mathcal{U}_{ad}^1\times\mathcal{U}_{ad}^2$ is called the set of admissible controls for the two players. We also introduce the following assumption:

\textbf{H1.} Functions $b,\sigma,\bar{\sigma}$ are continuously differentiable in $(x,x_\delta,v_1,v_2)$, $f$ is continuously differentiable in $(x,y,z,\bar{z},y_{\delta^+},v_1,v_2)$, $G$ is continuously differentiable in $x$. The all the partial derivatives of $b,\sigma,\bar{\sigma},f,G$ are uniformly bounded. %$b_x,b_{x_\delta},b_{v_i}$,$\sigma_x,\sigma_{x_\delta},\sigma_{v_i}$,$\bar{\sigma}_x,\bar{\sigma}_{x_\delta},\bar{\sigma}_{v_i}$,
%$f_x,f_y,f_z,f_{\bar{z}}$,
%$f_{y_{\delta^+}},f_{v_i}$,$(i=1,2)$ are uniformly bounded.
%Further, we assume there is constant $C$ such that $|h(t,x,v_1,v_2)|$$+|\sigma_1(t,x,v_1,v_2)|+|\sigma_2(t,x,v_1,v_2)|\leq C$ for $\forall (t,x,v_1,v_2)\in[0,T]\times\mathbb{R}\times U_1\times U_2$.

Then we have the following existence and uniqueness result which can be found in \cite{ChenWu10,PengYang09}.

\begin{theorem}
If $v_1(\cdot)$ and $v_2(\cdot)$ are admissible controls and assumption H1 holds, the AFBSDDE (\ref{afbsdde}) admits a unique solution $(x(\cdot),y(\cdot),z(\cdot),\bar{z}(\cdot))\in\mathbb{L}_{\mathcal{F}}^2(-\delta,T;\mathbb{R}^n)\times\mathbb{L}_{\mathcal{F}}^2(0,T+\delta;\mathbb{R}^m)
\times\mathbb{L}_{\mathcal{F}}^2(0,T;\mathbb{R}^{m\times d})\times\mathbb{L}_{\mathcal{F}}^2(0,T;\mathbb{R}^{m\times \bar{d}})$
\end{theorem}

The players have their own preferences which are described as the following cost functionals

\begin{equation*}
\begin{aligned}
J_i(v_1(\cdot),v_2(\cdot))=\mathbb{E}[\int_0^Tl_i(t,x^v(t),y^v(t),z^v(t),\bar{z}^v(t),v_{1}(t),v_{2}(t))dt+\Phi_i(x^v(T))+\gamma_i(y^v(0))].\\
\end{aligned}
\end{equation*}
Here $l_i: \Omega\times[0,T]\times\mathbb{R}^n\times\mathbb{R}^m\times\mathbb{R}^{m\times d}\times\mathbb{R}^{m\times \bar{d}}\times\mathbb{R}^{k_1}\times\mathbb{R}^{k_2}\rightarrow\mathbb{R}$, $\Phi_i: \Omega\times\mathbb{R}^n\rightarrow\mathbb{R}$, $\gamma_i: \Omega\times\mathbb{R}^m\rightarrow\mathbb{R}$$(i=1,2)$ are given continuous maps.
$l_i$, $\Phi_i$, and $\gamma_i$ satisfy the following condition:

\textbf{H2.}  Functions $l_i$, $\Phi_i$, and $\gamma_i$ are continuously differentiable with respect to $(x,y,z,\bar{z},v_1,v_2)$, $x$, and $y$ respectively. Moreover, there exists positive constant $C$ such that the partial derivatives of $l_i$, $\Phi_i$, and $\gamma_i$ are bounded by $C(1+|x|+|y|+|z|+|\bar{z}|+|v_1|+|v_2|)$, $C(1+|x|)$ and $C(1+|y|)$ respectively.
%$l_{ix},l_{iy},l_{iz},l_{i\bar{z}}$,$l_{iv_1},l_{iv_2},(i=1,2)$ are bounded by $C(1+|x|+|y|+|z|+|\bar{z}|+|v_1|+|v_2|)$.

Now we suppose that each player hopes to maximize his cost functional $J_i(v_1(\cdot),v_2(\cdot))$ by selecting a suitable admissible control $v_i(\cdot)(i=1,2)$. The problem is to find an admissible control $(u_1(\cdot),u_2(\cdot))\in\mathcal{U}_{ad}$ such that

\begin{equation}\label{J1J2}
\left\{
\begin{aligned}
J_1(u_1(\cdot),u_2(\cdot))=\sup\limits_{v_1(\cdot)\in \mathcal{U}_1}J_1(v_1(\cdot),u_2(\cdot)),\\
J_2(u_1(\cdot),u_2(\cdot))=\sup\limits_{v_2(\cdot)\in \mathcal{U}_2}J_2(u_1(\cdot),v_2(\cdot)).\\
\end{aligned}
\right.
\end{equation}

If we can find an admissible control $(u_1(\cdot),u_2(\cdot))$ satisfying (\ref{J1J2}), then we call it a Nash equilibrium point. In what follows, we aim to establish the necessary and sufficient condition for Nash equilibrium point subject to this game problem.

\section{Maximum principle}\label{MP}
In this section, we will establish a necessary condition (maximum principle) and a sufficient condition (verification theorem) for problem (\ref{J1J2}).

Let $(u_1(\cdot),u_2(\cdot))$ be an equilibrium point of the game problem. Then for any $0\leq\epsilon\leq 1$ and $(v_1(\cdot),v_2(\cdot))\in\mathcal{U}_{ad}$, we take the variational control $u_1^{\epsilon}(\cdot)=u_1(\cdot)+\epsilon v_1(\cdot)$ and $u_2^{\epsilon}(\cdot)=u_2(\cdot)+\epsilon v_2(\cdot)$. Because both $U_1$ and $U_2$ are convex, $(u_1^{\epsilon}(\cdot),u_2^{\epsilon}(\cdot)) $ is also in $\mathcal{U}_{ad}$. For simplicity, we denote by $(x^{u_1^\epsilon}(\cdot),y^{u_1^\epsilon}(\cdot),z^{u_1^\epsilon}(\cdot),\bar{z}^{u_1^\epsilon}(\cdot))$, $(x^{u_2^\epsilon}(\cdot),y^{u_2^\epsilon}(\cdot),z^{u_2^\epsilon}(\cdot),\bar{z}^{u_2^\epsilon}(\cdot))$, and $(x(\cdot),y(\cdot),z(\cdot),\bar{z}(\cdot))$ the corresponding state trajectories of system (\ref{afbsdde}) with control $(u_1^{\epsilon}(\cdot),u_2(\cdot))$, $(u_1(\cdot),u_2^{\epsilon}(\cdot))$ and $(u_1(\cdot),u_2(\cdot))$.

The following lemma gives an estimation of $(x(\cdot),y(\cdot),z(\cdot),\bar{z}(\cdot))$.

\begin{lem}\label{lem} Let $H1$ hold. For $i=1,2$,
\begin{equation*}
\sup\limits_{0\leq t\leq T}\mathbb{E}|x^{u_i^\epsilon}(t)-x^{}(t)|^2\leq C\epsilon^2,
\end{equation*}

\begin{equation*}
\sup\limits_{0\leq t\leq T}\mathbb{E}|y^{u_i^\epsilon}(t)-y^{}(t)|^2\leq C\epsilon^2,
\end{equation*}

\begin{equation*}
\mathbb{E}\int_0^T|z^{u_i^\epsilon}(t)-z^{}(t)|^2dt\leq C\epsilon^2,
\end{equation*}

\begin{equation*}
\mathbb{E}\int_0^T|\bar{z}^{u_i^\epsilon}(t)-\bar{z}^{}(t)|^2dt\leq C\epsilon^2,
\end{equation*}
\end{lem}

\textbf{Proof.} Using It\^{o}'s formula to $|x^{u_i^\epsilon}(t)-x^{}(t)|^2$ and Gronwall's inequality, we draw the conclusion.

For notation simplicity, we set $\zeta(t)=\zeta(t,x(t),x_\delta(t),u_{1}(t),u_{2}(t))\ \text{for}\  \zeta=b, \sigma, \bar{\sigma}$; $f(t)=f(t,x(t),y(t)$,
$z(t),\bar{z}(t),y_{\delta^+}(t),u_{1}(t),u_{2}(t))$, and $l_i(t)=l_i(t,x(t),y(t),z(t),\bar{z}(t),u_{1}(t),u_{2}(t))(i=1,2)$.

We introduce the following variational equations:
\begin{equation}\label{Variational xy}
\left\{
\begin{aligned}
dx_i^1(t)=\ &[b_x(t)x_i^1(t)+b_{x_\delta}(t)x_i^1(t-\delta)+b_{v_i}(t)v_i(t)]dt+[\sigma_x(t)x_i^1(t)+\sigma_{x_\delta}(t)x_i^1(t-\delta)\\
&+\sigma_{v_i}(t)v_i(t)]dW(t)+[\bar{\sigma}_x(t)x_i^1(t)+\bar{\sigma}_{x_\delta}(t)x_i^1(t-\delta)+\bar{\sigma}_{v_i}(t)v_i(t)]d\bar{W}(t),\\
-dy_i^1(t)=\ &\{f_x(t)x_i^1(t)+f_y(t)y_i^1(t)+f_z(t)z_i^1(t)+f_{\bar{z}}(t)\bar{z}_i^1(t)+\mathbb{E}^{\mathcal{F}_t}[f_{y_{\delta^+}}(t)y_{i}^1(t+\delta)]\\
&+f_{v_i}(t)v_i(t)\}dt-z_i^1(t)dW(t)-\bar{z}_i^1(t)d\bar{W}(t),\quad t\in[0,T],\\
x_i^1(t)=\ &0,\quad t\in[-\delta,0],\\
y_i^1(T)=\ &G_x(x(T))x_i^1(T),\ y_i^1(t)=0,\quad t\in(T,T+\delta],\quad(i=1,2).
\end{aligned}
\right.
\end{equation}

Next, setting

\[
\phi_i^{\epsilon}(t)=\frac{\phi^{u_i^{\epsilon}}(t)-\phi(t)}{\epsilon}-\phi_i^1(t),\ \ \text{for}\ \ \phi=x,y,z,\bar{z},\quad(i=1,2),
\]

Then we can get the following two lemmas by using Lemma \ref{lem}. The technique is classical (see Chen and Wu \cite{ChenWu10}). Thus we omit the details and only state the main result for simplicity.

\begin{lem}
Let H1 hold. For $i=1,2$,
\begin{equation*}
\lim\limits_{\epsilon\rightarrow0}\sup\limits_{0\leq t\leq T}\mathbb{E}|x_i^\epsilon(t)|^2=0,
\end{equation*}

\begin{equation*}
\lim\limits_{\epsilon\rightarrow0}\sup\limits_{0\leq t\leq T}\mathbb{E}|y_i^\epsilon(t)|^2=0,
\end{equation*}

\begin{equation*}
\lim\limits_{\epsilon\rightarrow0}\mathbb{E}\int_0^T|z_i^\epsilon(t)|^2dt=0,
\end{equation*}

\begin{equation*}
\lim\limits_{\epsilon\rightarrow0}\mathbb{E}\int_0^T|\bar{z}_i^\epsilon(t)|^2dt=0.
\end{equation*}
\end{lem}

\begin{lem}
Let H1 and H2 hold. For $i=1,2$,
\begin{equation}\label{var equality}
\begin{aligned}
\mathbb{E}\int_0^T&[l_{ix}^\tau (t) x_i^1(t)+l_{iy}^\tau (t) y_i^1(t)+l_{iz}^\tau (t) z_i^1(t)+l_{i\bar{z}}^\tau (t) \bar{z}_i^1(t)+l_{iv_i}^\tau (t) v_i^1(t)]dt\\
&+\mathbb{E}[\Phi_{ix}^\tau (x(T))x_i^1(T)]+\gamma_{iy}^\tau (y(0))y_i^1(0)\leq 0.
\end{aligned}
\end{equation}
\end{lem}

We introduce the adjoint equation as

\begin{equation}\label{adjoint}
\left\{
\begin{aligned}
dp_i(t)=\ &[f_y^\tau(t)p_i(t)+f_{y_{\delta^+}}^\tau(t-\delta)p_i(t-\delta)-l_{iy}(t)]dt+[f_z^\tau(t)p_i(t)-l_{iz}(t)]dW(t)\\
&+[f_{\bar{z}}^\tau(t)p_i(t)-l_{i\bar{z}}(t)]d\bar{W}(t),\\
-dq_i(t)=\ &\{b_x^\tau(t)q_i(t)+\sigma_x^\tau(t)k_i(t)+\bar{\sigma}_x^\tau(t)\bar{k}_i(t)-f_x^\tau(t)p_i(t)+\mathbb{E}^{\mathcal{F}_t}[b_{x_\delta}^\tau(t+\delta)q_i(t+\delta)\\
&+\sigma_{x_\delta}^\tau(t+\delta)k_i(t+\delta)+\bar{\sigma}_{x_\delta}^\tau(t+\delta)\bar{k}_i(t+\delta)]+l_{ix}(t)\}dt-k_i(t)dW(t)-\bar{k}_i(t)d\bar{W}(t),\\
p_i(0)=&-\gamma_y(y(0)),\ p_i(t)=0,\quad t\in[-\delta,0),\\
q_i(T)=&-G_x^\tau(x(T))p_i(T)+\Phi_{ix}(x(T)),\ q_i(t)=k_i(t)=\bar{k}_i(t)=0,\quad t\in(T,T+\delta],\quad (i=1,2).
\end{aligned}
\right.
\end{equation}

This equation is also an AFBSDDE. By the existence and uniqueness result in \cite{ChenWu10,PengYang09}, we know that (\ref{adjoint}) admits a unique solution $(p_i(t),q_i(t),k_i(t),\bar{k}_i(t))(i=1,2)$.

Define the Hamiltonian function $H_i$ by

\begin{equation*}
\begin{aligned}
H_i&(t,x,y,z,\bar{z},x_\delta,y_{\delta^+},v_1,v_2;p_i,q_i,k_i,\bar{k}_i)=\langle q_i,b(t,x,x_\delta,v_1,v_2)\rangle+\langle k_i,\sigma(t,x,x_\delta,v_1,v_2)\rangle\\
+&\langle \bar{k}_i,\bar{\sigma}(t,x,x_\delta,v_1,v_2)\rangle-\langle p_i,f(t,x,y,z,\bar{z},y_{\delta^+},v_1,v_2)\rangle+l_i(t,x,y,z,\bar{z},v_1,v_2),\quad (i=1,2).
\end{aligned}
\end{equation*}

Then (\ref{adjoint}) can be rewritten as the following stochastic Hamiltonian system's type:

\begin{equation*}
\left\{
\begin{aligned}
dp_i(t)=\ &[-H_{iy}(t)-H_{iy_{\delta^+}}(t-\delta)]dt-H_{iz}(t)dW(t)-H_{i\bar{z}}(t)d\bar{W}(t),\\
-dq_i(t)=\ &\{H_{ix}(t)+\mathbb{E}^{\mathcal{F}_t}[H_{ix_\delta}(t+\delta)]\}dt-k_i(t)dW(t)-\bar{k}_i(t)d\bar{W}(t),\quad t\in[0,T],\\
p_i(0)=\ &-\gamma_y(y(0)),\ p_i(t)=0,\quad t\in[-\delta,0),\\
q_i(T)=\ &-G_x^\tau(x(T))p_i(T)+\Phi_{ix}(x(T)),\ q_i(t)=k_i(t)=\bar{k}_i(t)=0,\quad t\in(T,T+\delta],\quad (i=1,2).
\end{aligned}
\right.
\end{equation*}
where $H_i(t)=H_i(t,x(t),y(t),z(t),\bar{z}(t),x_\delta(t),y_{\delta^+}(t),v_1(t),v_2(t);p_i(t),q_i(t),k_i(t),\bar{k}_i(t))$.

\begin{theorem}\label{nec MP}
Let H1 and H2 hold. Suppose that $(u_1(\cdot),u_2(\cdot))$ is an equilibrium point of our problem and $(x(\cdot),y(\cdot),z(\cdot),\bar{z}(\cdot))$ is the corresponding state trajectory. Then we have

\begin{equation*}
\mathbb{E}[\langle H_{iv_i}(t),v_i-u_i(t)\rangle|\mathcal{G}_t^i]\leq 0,\quad(i=1,2)
\end{equation*}
for any $v_i\in U_i$ a.e., where $(p_i(\cdot),q_i(\cdot),k_i(\cdot),\bar{k}_i(\cdot)),(i=1,2)$ is the solution of the adjoint equation (\ref{adjoint}).

\textbf{Proof}.
Apply It\^{o}'s formula to $\langle q_1(\cdot),x_1^1(\cdot)\rangle$, we get

\begin{equation}\label{q1x1}
\begin{aligned}
&\ \ \ \ \mathbb{E}\langle-G_x^\tau(x(T))p_1(T)+\Phi_{1x}(x(T)),x_1^1(T)\rangle\\
&=\mathbb{E}\int_0^T[\langle f_x^\tau(t)p_1(t),x_1^1(t)\rangle+\langle b_{x_\delta}^\tau(t)q_1(t),x_1^1(t-\delta)\rangle-\langle \mathbb{E}^{\mathcal{F}_t}[b_{x_\delta}^\tau(t+\delta)q_1(t+\delta)],x_1^1(t)\rangle\\
&\ \ \ +\langle \sigma_{x_\delta}^\tau(t)k_1(t)+\bar{\sigma}_{x_\delta}^\tau(t)\bar{k}_1(t),x_1^1(t-\delta)\rangle-\langle \mathbb{E}^{\mathcal{F}_t}[\sigma_{x_\delta}^\tau(t+\delta)k_1(t+\delta)+\bar{\sigma}_{x_\delta}^\tau(t+\delta)\bar{k}_1(t+\delta)],x_1^1(t)\rangle\\
&\ \ \ +\langle q_1(t),b_{v_1}(t)v_1(t)\rangle+\langle k_1(t),\sigma_{v_1}(t)v_1(t)\rangle+\langle \bar{k}_1(t),\bar{\sigma}_{v_1}(t)v_1(t)\rangle-\langle l_{1x}(t),x_1^1(t)\rangle]dt.
\end{aligned}
\end{equation}

Noticing the initial and terminal conditions, we have

\begin{equation*}
  \begin{aligned}
    &\ \ \ \ \mathbb{E}\int_0^T[\langle b_{x_\delta}^\tau(t)q_1(t),x_1^1(t-\delta)\rangle-\langle \mathbb{E}^{\mathcal{F}_t}[b_{x_\delta}^\tau(t+\delta)q_1(t+\delta)],x_1^1(t)\rangle] dt\\
    &=\mathbb{E}\int_0^T\langle b_{x_\delta}^\tau(t)q_1(t),x_1^1(t-\delta)\rangle dt-\mathbb{E}\int_\delta^{T+\delta}\langle b_{x_\delta}^\tau(t)q_1(t),x_1^1(t-\delta)\rangle dt\\
    &=\mathbb{E}\int_0^\delta\langle b_{x_\delta}^\tau(t)q_1(t),x_1^1(t-\delta)\rangle dt-\mathbb{E}\int_T^{T+\delta}\langle b_{x_\delta}^\tau(t)q_1(t),x_1^1(t-\delta)\rangle dt\\
    &=0.
  \end{aligned}
\end{equation*}

Similarly, we also have

\begin{equation*}
  \begin{aligned}
  \mathbb{E}\int_0^T[\langle \sigma_{x_\delta}^\tau(t)k_1(t)+\bar{\sigma}_{x_\delta}^\tau(t)\bar{k}_1(t),x_1^1(t-\delta)\rangle-\langle \mathbb{E}^{\mathcal{F}_t}[\sigma_{x_\delta}^\tau(t+\delta)k_1(t+\delta)+\bar{\sigma}_{x_\delta}^\tau(t+\delta)\bar{k}_1(t+\delta)],x_1^1(t)\rangle] dt=0.
  \end{aligned}
\end{equation*}

Apply It\^{o}'s formula to $\langle p_1(\cdot),y_1^1(\cdot)\rangle$,

\begin{equation}\label{p1y1}
\begin{aligned}
&\ \ \ \ \mathbb{E}\langle p_1(T),G_x(x(T))x_1^1(T)\rangle+\langle \gamma_y(y(0)),y_1^1(0)\rangle\\
&=\mathbb{E}\int_0^T[\langle f_{y_{\delta^+}}^\tau(t-\delta)p_1(t-\delta),y_1^1(t)\rangle-\langle p_1(t),\mathbb{E}^{\mathcal{F}_t}[f_{y_{\delta^+}}(t)y_1^1(t+\delta)]\rangle\\
&\ \ \ -\langle p_1(t),f_x(t)x_1^1(t)+f_{v_1}(t)v_1(t)\rangle-\langle l_{1y}(t),y_1^1(t)\rangle-\langle l_{1z}(t),z_1^1(t)\rangle-\langle l_{1\bar{z}}(t),\bar{z}_1^1(t)\rangle]dt.
\end{aligned}
\end{equation}

Noticing the initial and terminal conditions, we have

\begin{equation*}
  \begin{aligned}
    &\ \ \ \ \mathbb{E}\int_0^T[\langle f_{y_{\delta^+}}^\tau(t-\delta)p_1(t-\delta),y_1^1(t)\rangle-\langle p_1(t),\mathbb{E}^{\mathcal{F}_t}[f_{y_{\delta^+}}(t)y_1^1(t+\delta)]\rangle] dt\\
    &=\mathbb{E}\int_0^T\langle f_{y_{\delta^+}}^\tau(t-\delta)p_1(t-\delta),y_1^1(t)\rangle dt-\mathbb{E}\int_\delta^{T+\delta}\langle f_{y_{\delta^+}}^\tau(t-\delta)p_1(t-\delta),y_1^1(t)\rangle dt\\
    &=\mathbb{E}\int_0^\delta\langle f_{y_{\delta^+}}^\tau(t-\delta)p_1(t-\delta),y_1^1(t)\rangle dt-\mathbb{E}\int_T^{T+\delta}\langle f_{y_{\delta^+}}^\tau(t-\delta)p_1(t-\delta),y_1^1(t)\rangle dt\\
    &=0.
  \end{aligned}
\end{equation*}

From (\ref{q1x1}) and (\ref{p1y1}), we have

\begin{equation}\label{qyqx}
\begin{aligned}
&\ \ \ \ \mathbb{E}\langle\Phi_{1x}(x(T)),x_1^1(T)\rangle+\langle \gamma_y(y(0)),y_1^1(0)\rangle\\
&=\mathbb{E}\int_0^T[\langle q_1(t),b_{v_1}(t)v_1(t)\rangle+\langle k_1(t),\sigma_{v_1}(t)v_1(t)\rangle+\langle \bar{k}_1(t),\bar{\sigma}_{v_1}(t)v_1(t)\rangle-\langle p_1(t),f_{v_1}(t)v_1(t)\rangle\\
&\ \ \ -\langle l_{1y}(t),y_1^1(t)\rangle-\langle l_{1z}(t),z_1^1(t)\rangle-\langle l_{1\bar{z}}(t),\bar{z}_1^1(t)\rangle-\langle l_{1x}(t),x_1^1(t)\rangle]dt.
\end{aligned}
\end{equation}

Substituting (\ref{qyqx}) into (\ref{var equality}), it follows that

\begin{equation*}
\mathbb{E}\int_0^T\langle H_{1v_1}(t),v_1(t)\rangle dt\leq 0
\end{equation*}
for any $v_1(\cdot)$ such that $u_1(\cdot)+v_1(\cdot)\in\mathcal{U}_{ad}^1$. If we let $\nu_1(\cdot)=u_1(\cdot)+v_1(\cdot)$, then above equation implies that

\begin{equation*}
\mathbb{E}\langle H_{1v_1}(t),\nu_1(t)-u_1(t)\rangle \leq 0.
\end{equation*}

Furthermore, we set

\begin{equation*}
\omega_1(t)=v_1 1_A+u_1(t)1_{\Omega-A},\quad \forall v_1\in U_1,\quad\forall A\in\mathcal{G}_t^1,
\end{equation*}
then it is obvious that $\omega_1(\cdot)\in\mathcal{U}_{ad}^1$ and $v_1(t)=(v_1-u_1(t))1_A$. So
\begin{equation*}
\mathbb{E}[1_A\langle H_{1v_1}(t),v_1-u_1(t)\rangle]\leq 0
\end{equation*}
for any $A\in\mathcal{G}_t^1.$ This implies

\begin{equation*}
\mathbb{E}[\langle H_{1v_1}(t),v_1-u_1(t)\rangle|\mathcal{G}_t^1]\leq 0,\quad a.e.
\end{equation*}
for any $v_1\in U_1$.

Repeating the same process to deal with the case $i=2$, we can show that the other equality also holds for any $v_2\in U_2$. Our proof is completed.

\end{theorem}

\begin{remark}
If $(u_1(\cdot),u_2(\cdot))$ is an equilibrium point of non-zero sum differential game and $(u_1(\cdot),u_2(\cdot))$ is an interior point of $U_1\times U_2$ for all $t\in[0,T]$, then the inequality in Theorem \ref{nec MP} are equivalent to the following equations:
\begin{equation*}
\mathbb{E}[\langle H_{iv_i}(t),v_i-u_i(t)\rangle|\mathcal{G}_t^i]=0, \quad \forall v_i\in U_i\quad a.e., \quad (i=1,2).
\end{equation*}
\end{remark}

On the other hand, we will aim to build a sufficient maximum principle called verification theorem for equilibrium point under some concavity assumptions of $H_i$. At this moment, assumption H2 can be relaxed to

\textbf{H3.}  Functions $l_i$, $\Phi_i$, and $\gamma_i$ are differentiable with respect to $(x,y,z,\bar{z},v_1,v_2)$, $x$, and $y$ respectively satisfying the condition that for each $(v_1(
\cdot),v_2(\cdot))\in\mathcal{U}_{ad}$, $l_i(\cdot,x^v(t),y^v(t),z^v(t),\bar{z}^v(t),v_{1}(t),v_{2}(t))\in\mathbb{L}^1_{\mathcal{F}}(0,T;\mathbb{R})$.

\begin{theorem}\label{suf MP}
Let H1 and H3 hold. Let $(u_1(\cdot),u_2(\cdot))\in\mathcal{U}_{ad}^1\times\mathcal{U}_{ad}^2$ be given and $(x(\cdot),y(\cdot),z(\cdot),\bar{z}(\cdot))$ be the corresponding trajectory.

Suppose
\begin{equation*}
\begin{aligned}
\mathbb{E}[{H}_1(t)|\mathcal{G}_t^1]=\sup\limits_{v_1\in U_1}\mathbb{E}[{H}_1^{v_1}(t)|\mathcal{G}_t^1],\\
\mathbb{E}[{H}_2(t)|\mathcal{G}_t^2]=\sup\limits_{v_2\in U_2}\mathbb{E}[{H}_2^{v_2}(t)|\mathcal{G}_t^2],
\end{aligned}
\end{equation*}

where
\begin{equation*}
\begin{aligned}
{H}_{1}^{v_1}(t)=& H_{1}(t,x(t),y(t),z(t),\bar{z}(t),x_\delta(t),y_{\delta^+}(t),v_1(t),u_2(t);p_1(t),q_1(t),k_1(t),\bar{k}_{1}(t)),\\
{H}_{2}^{v_2}(t)=& H_{2}(t,x(t),y(t),z(t),\bar{z}(t),x_\delta(t),y_{\delta^+}(t),u_1(t),v_2(t);p_2(t),q_2(t),k_2(t),\bar{k}_{2}(t)).
\end{aligned}
\end{equation*}

Suppose $\mathbb{E}[{H}_{iv_i}^{v_i}(t)|\mathcal{G}_t^i]$ is continuous at $v_i=u_i(t)(i=1,2)$ for any $t\in[0,T]$.

Suppose
\begin{equation*}
\begin{aligned}
(x,y,z,\bar{z},x_\delta,y_{\delta^+},v_i)\mapsto& H_i^{v_i}(t)\quad(i=1,2),\\
x\mapsto& \Phi_i(x)\ \quad(i=1,2),\\
y\mapsto& \gamma_i(y)\ \ \quad(i=1,2)
\end{aligned}
\end{equation*}
are concave functions respectively, and $G(x)=M_Tx, M_T\in\mathbb{R}^{m\times n},\forall x\in\mathbb{R}^n$. Then $(u_1(\cdot),u_2(\cdot))$ is an equilibrium point.
\end{theorem}

\textbf{Proof}. For any $v_1(\cdot)\in\mathcal{U}_{ad}^1$, let $(x^{v_1}(\cdot),y^{v_1}(\cdot),z^{v_1}(\cdot),\bar{z}^{v_1}(\cdot))$ be the trajectory corresponding to the control $(v_1(\cdot),u_2(\cdot))\in\mathcal{U}_{ad}$. We consider

$$
J_1(v_1(\cdot),u_2(\cdot))-J_1(u_1(\cdot),u_2(\cdot))=A+B+C,
$$
with

\begin{equation*}
\begin{aligned}
A=&\mathbb{E}\int_0^T[l_1(t,\Theta^{v_1}(t),v_1(t),u_2(t))-l_1(t,\Theta(t),u_1(t),u_2(t))]dt,\\
B=&\mathbb{E}[\Phi_1(x^{v_1}(T))-\Phi_1(x(T))],\\
C=&\gamma_1(y^{v_1}(0))-\gamma_1(y(0)).
\end{aligned}
\end{equation*}
where $\Theta(t)=(x(t),y(t),z(t),\bar{z}(t))$ and $\Theta^{v_1}(t)=(x^{v_1}(t),y^{v_1}(t),z^{v_1}(t),\bar{z}^{v_1}(t))$.

Since $\gamma_1$ is concave on $y$, then

$$
C\leq\gamma_{1y}^\tau(y(0))(y^{v_1}(0)-y(0)).
$$

Applying It\^{o}'s formula to $\langle p_1(\cdot),y^{v_1}(\cdot)-y(\cdot)\rangle$ and taking expectation, we get

\begin{equation}\label{C}
  \begin{aligned}
    &C\leq\mathbb{E}\int_0^T[-\langle p_1(t),f^{v_1}(t)-f(t)\rangle-\langle H_{1y}(t)+H_{1y_{\delta^+}}(t-\delta),y^{v_1}(t)-y(t)\rangle\\
    &-\langle H_{1z}(t),z^{v_1}(t)-z(t)\rangle-\langle H_{1\bar{z}}(t),\bar{z}^{v_1}(t)-\bar{z}(t)\rangle]dt-\mathbb{E}\langle p_1(T),M_T(x^{v_1}(T)-x(T))\rangle
  \end{aligned}
\end{equation}
where $f(t)=f(t,\Theta(t),\mathbb{E}^{\mathcal{F}_t}[y_{\delta^+}(t)],u_1(t),u_2(t))$ and $f^{v_1}(t)=f(t,\Theta^{v_1}(t),\mathbb{E}^{\mathcal{F}_t}[y_{\delta^+}^{v_1}(t)],v_1(t),u_2(t))$.

Due to $\Phi_1$ is concave on $x$, then

$$
B\leq\mathbb{E}\Phi_{1x}^\tau(x(T))(x^{v_1}(T)-x(T)).
$$

Applying It\^{o}'s formula to $\langle q_1(\cdot),x^{v_1}(\cdot)-x(\cdot)\rangle$ and taking expectation, we get

\begin{equation}\label{B}
  \begin{aligned}
    &B\leq\mathbb{E}\int_0^T[\langle q_1(t),b^{v_1}(t)-b(t)\rangle+\langle k_1(t),\sigma^{v_1}(t)-\sigma(t)\rangle+\langle \bar{k}_1(t),\bar{\sigma}^{v_1}(t)-\bar{\sigma}(t)\rangle\\
    &-\langle H_{1x}(t)+\mathbb{E}^{\mathcal{F}_t}[H_{1x_{\delta}}(t+\delta)],x^{v_1}(t)-x(t)\rangle ]dt+\mathbb{E}\langle M_T^\tau p_1(T),x^{v_1}(T)-x(T)\rangle
  \end{aligned}
\end{equation}
where $b(t)=b(t,x(t),x_\delta(t),u_1(t),u_2(t))$ and $b^{v_1}(t)=b(t,x^{v_1}(t),x^{v_1}_\delta(t),v_1(t),u_2(t))$, etc.

Moreover, we have

\begin{equation}\label{A}
  \begin{aligned}
    A=\ &\mathbb{E}\int_0^T[H_1^{v_1}(t)-H_1(t)]dt-\mathbb{E}\int_0^T[\langle q_1(t),b^{v_1}(t)-b(t)\rangle+\langle k_1(t),\sigma^{v_1}(t)-\sigma(t)\rangle\\
    &+\langle \bar{k}_1(t),\bar{\sigma}^{v_1}(t)-\bar{\sigma}(t)\rangle-\langle p_1(t),f^{v_1}(t)-f(t)\rangle]dt.
  \end{aligned}
\end{equation}

From (\ref{C})-(\ref{A}), we can obtain
$$
\begin{aligned}
J&_1(v_1(\cdot),u_2(\cdot))-J_1(u_1(\cdot),u_2(\cdot))=A+B+C\\
&\leq \mathbb{E}\int_0^T[(H_1^{v_1}(t)-H_1(t))-\langle H_{1x}(t)+\mathbb{E}^{\mathcal{F}_t}[H_{1x_{\delta}}(t+\delta)],x^{v_1}(t)-x(t)\rangle\\
&-\langle H_{1y}(t)+H_{1y_{\delta^+}}(t-\delta),y^{v_1}(t)-y(t)\rangle-\langle H_{1z}(t),z^{v_1}(t)-z(t)\rangle-\langle H_{1\bar{z}}(t),\bar{z}^{v_1}(t)-\bar{z}(t)\rangle]dt.
\end{aligned}
$$

Note that

$$
\begin{aligned}
&\ \ \ \ \mathbb{E}\int_0^T\langle H_{1x_\delta}(t),x^{v_1}(t-\delta)-x(t-\delta)\rangle dt-\mathbb{E}\int_0^T\langle\mathbb{E}^{\mathcal{F}_t}[H_{1x_\delta}(t+\delta),x^{v_1}(t)-x(t)\rangle dt\\
&=\mathbb{E}\int_{-\delta}^0\langle H_{1x_\delta}(t+\delta),x^{v_1}(t)-x(t)\rangle dt-\mathbb{E}\int_{T-\delta}^T\langle H_{1x_\delta}(t+\delta),x^{v_1}(t)-x(t)\rangle dt\\
&=0,
\end{aligned}
$$
due to the fact that $x^{v_1}(t)=x(t)=\xi$ for any $t\in[-\delta,0)$ and $H_{1x_\delta}(t)=0$ for any $t\in(T,T+\delta]$.

Similarly, we have

$$
\begin{aligned}
&\ \ \ \ \mathbb{E}\int_0^T\langle H_{1y_{\delta^+}}(t),\mathbb{E}^{\mathcal{F}_t}[y^{v_1}(t+\delta)-y(t+\delta)]\rangle dt-\mathbb{E}\int_0^T\langle H_{1y_{\delta^+}}(t-\delta),y^{v_1}(t)-y(t)\rangle dt\\
&=\mathbb{E}\int_T^{T+\delta}\langle H_{1y_{\delta^+}}(t-\delta),y^{v_1}(t)-y(t)\rangle dt-\mathbb{E}\int_0^\delta\langle H_{1y_{\delta^+}}(t-\delta),y^{v_1}(t)-y(t)\rangle dt\\
&=0,
\end{aligned}
$$
due to the fact that $y^{v_1}(t)=y(t)=\varphi(t)$ for any $t\in(T,T+\delta]$ and $H_{1y_{\delta^+}}(t)=0$ for any $t\in[-\delta,0)$.

By the concavity of $H_1$, we derive that

$$
\begin{aligned}
J_1(v_1(\cdot),u_2(\cdot))-J_1(u_1(\cdot),u_2(\cdot))&\leq\mathbb{E}\int_0^T\langle H_{1v_1}(t),v_1(t)-u_1(t)\rangle dt\\
&=\mathbb{E}\int_0^T\mathbb{E}[\langle H_{1v_1}(t),v_1(t)-u_1(t)\rangle|\mathcal{G}_t^1]dt.
\end{aligned}
$$

Because for any $t\in[0,T]$, $v_1\mapsto\mathbb{E}[H_1^{v_1}(t)|\mathcal{G}_t^1]$ is maximal at $v_1=u(t)$ and $H_{1v_1}^{v_1}(t)$ is continuous at $v_1=u_1(t)$, then we have

$$
\mathbb{E}[\langle H_{1v_1}(t),v_1(t)-u_1(t)\rangle|\mathcal{G}_t^1]\leq 0,\quad a.e.
$$

It follows that
$$
J_1(u_1(\cdot),u_2(\cdot))=\sup\limits_{v_1\in\mathcal{U}_1}J_1(v_1(\cdot),u_2(\cdot)).
$$

Repeating the same process to deal with case $i=2$, we can obtain
$$
J_2(u_1(\cdot),u_2(\cdot))=\sup\limits_{v_2\in\mathcal{U}_2}J_2(u_1(\cdot),v_2(\cdot)).
$$

Hence we draw the desired conclusion.

In conclusion, with the help of Theorem \ref{nec MP} and Theorem \ref{suf MP}, we can formally solve the Nash equilibrium point $(u_1(\cdot),u_2(\cdot))$. We can first use the necessary principle to get the candidate equilibrium point and then use the verification theorem to check whether the candidate point is the equilibrium one. Let us discuss a linear-quadratic case.

\section{A linear quadratic case}
In this section, we study a linear-quadratic case, which can be seen as a special case of the general system discussed in Section \ref{MP} and aim to give a unique Nash equilibrium point explicitly. For notational simplification, we assume the dimension of Brownian motion $d=\bar{d}=1$ and notations are the same as the above sections if there is no specific illustration.

Consider a linear game system with delayed and anticipated states:
\begin{equation}\label{afbsdde lq}
\left\{
\begin{aligned}
dx^v(t)=\ &[A(t)x^v(t)+\bar{A}(t)x^v_{\delta}(t)+B_{1}(t)v_{1}(t)+B_{2}(t)v_{2}(t)]dt+[C(t)x^v(t)+\bar{C}(t)x^v_{\delta}(t)\\
&+D_{1}(t)v_{1}(t)+D_{2}(t)v_{2}(t)]dW(t),\\
-dy^v(t)=\ &[E(t)x^v(t)+F(t)y^v(t)+G(t)z^v(t)+\bar{G}(t)\bar{z}^v(t)+\bar{F}(t)y^v_{\delta^+}(t)+H_{1}(t)v_{1}(t)\\
&+H_{2}(t)v_{2}(t)]dt-z^v(t)dW(t)-\bar{z}^v(t)d\bar{W}(t),\quad t\in[0,T],\\
x^v(t)=\ &\xi(t),\quad t\in[-\delta,0],\\
y^v(T)=\ &M_Tx^v(T),\ y^v(t)=\varphi(t),\quad t\in(T,T+\delta].\\
\end{aligned}
\right.
\end{equation}
where all the coefficients are bounded, deterministic matrices defining on $[0,T]$, $\xi(\cdot)\in\mathbb{C}([-\delta,0];\mathbb{R}^n)$, $\varphi(\cdot)\in\mathbb{L}^2_\mathcal{F}(T,T+\delta;\mathbb{R}^m)$. For any given $(v_1(\cdot),v_2(\cdot))\in\mathcal{U}_{ad}$, it is easy to show know that (\ref{afbsdde lq}) admits a unique solution $(x^v(\cdot),y^v(\cdot),z^v(\cdot),\bar{z}^v(\cdot))$. Here we only consider the case that $x^v(\cdot)$ is driven by one Brownian motion $W(\cdot)$ just for notation simplicity. All the technique and proof is similar.

In addition, two players aim to maximize their index functionals for $i=1,2$:

\begin{equation*}
\begin{aligned}
J_i(v_1(\cdot),v_2(\cdot))=\ &\frac{1}{2}\mathbb{E}[\int_0^T[\langle O_i(t)x^v(t),x^v(t)\rangle+\langle P_i(t)y^v(t),y^v(t)\rangle+\langle Q_i(t)z^v(t),z^v(t)\rangle\\
&+\langle \bar{Q}_i(t)\bar{z}^v(t),\bar{z}^v(t)\rangle+\langle R_i(t)v_i(t),v_i(t)\rangle]dt+\langle M_ix^v(t),x^v(t)\rangle+\langle N_iy^v(0),y^v(0)\rangle].\\
\end{aligned}
\end{equation*}
where $O_i(\cdot)$, $P_i(\cdot)$, $Q_i(\cdot)$, $\bar{Q}_i(\cdot)$ are bounded deterministic non-positive symmetric matrices, $R_i(\cdot)$ is  bounded deterministic negative symmetric matrices, $R_i^{-1}(\cdot)$ is bounded,  $M_i$, $N_i$ are deterministic non-positive symmetric matrices for $i=1,2$.

According to Theorem \ref{nec MP}, the Hamiltonian function is given by
\begin{equation*}
\begin{aligned}
\mathbb{H}&_i(t,x,y,z,\bar{z},x_\delta,y_{\delta^+},v_1,v_2;p_i,q_i,k_i)=\langle q_i,A(t)x+\bar{A}(t)x_\delta+B_{1}(t)v_{1}+B_{2}(t)v_{2}\rangle\\
&+\langle k_i,C(t)x+\bar{C}(t)x_\delta+D_{1}(t)v_{1}+D_{2}(t)v_{2}\rangle-\langle p_i,E(t)x+F(t)y+G(t)z+\bar{G}(t)\bar{z}+\bar{F}(t)y_{\delta^+}\\
&+H_{1}(t)v_{1}+H_{2}(t)v_{2}\rangle+\frac{1}{2}[\langle O_i(t)x,x\rangle+\langle P_i(t)y,y\rangle+\langle Q_i(t)z,z\rangle+\langle \bar{Q}_i(t)\bar{z},\bar{z}\rangle+\langle R_i(t)v_i,v_i\rangle].\\
\end{aligned}
\end{equation*}

If $(u_1(\cdot),u_2(\cdot))$ is the Nash equilibrium point, then

\begin{equation}\label{ui lq}
\begin{aligned}
u_i(t)=-R^{-1}_i(t)[B^{\tau}_i(t)\hat{q}_i(t)+D^{\tau}_i(t)\hat{k}_i(t)-H^{\tau}_i(t)\hat{p}_i(t)],\quad t\in[0,T],\quad (i=1,2),\\
\end{aligned}
\end{equation}
where $\hat{q}_i(t)=\mathbb{E}[q_i(t)|{\mathcal{G}_t}]$, etc., and $(p_i(\cdot),q_i(\cdot),k_i(\cdot))$ is the solution of the following adjoint equation:

\begin{equation}\label{ad lq}
\left\{
\begin{aligned}
dp_i(t)=\ &[F^\tau(t)p_i(t)+\bar{F}^\tau(t-\delta)p_i(t-\delta)-P_{i}(t)y(t)]dt+[G^\tau(t)p_i(t)-Q_{i}(t)z(t)]dW(t)\\
&+[\bar{G}^\tau(t)p_i(t)-\bar{Q}_{i}(t)\bar{z}(t)]d\bar{W}(t),\\
-dq_i(t)=\ &\{A^\tau(t)q_i(t)+C^\tau(t)k_i(t)-E^\tau(t)p_i(t)+\mathbb{E}^{\mathcal{F}_t}[\bar{A}^\tau(t+\delta)q_i(t+\delta)+\bar{C}^\tau(t+\delta)k_i(t+\delta)]\\
&+O_i(t)x(t)\}dt-k_i(t)dW(t),\quad t\in[0,T],\\
p_i(0)=&-N_iy(0),\ p_i(t)=0,\quad t\in[-\delta,0),\\
q_i(T)=&-M_Tp_i(T)+M_ix(T),\ q_i(t)=0,\quad t\in(T,T+\delta],\quad (i=1,2).
\end{aligned}
\right.
\end{equation}

We know that the setting $\mathcal{G}_t\subseteq\mathcal{F}_t$ is very general. In order to get an explicit expression of the equilibrium point, we suppose $\mathcal{G}_t=\sigma\{W(s);0\leq s\leq t\}$ in the rest of this section.

We denote the filtering of state process $x(t)$ by $\hat{x}(t)=\mathbb{E}[x(t)|{\mathcal{G}_t}]$, etc, and note that $\mathbb{E}[y(t+\delta)|\mathcal{G}_t]=\mathbb{E}\{[y(t+\delta)|\mathcal{G}_{t+\delta}]|\mathcal{G}_t\}=\mathbb{E}[\hat{y}(t+\delta)|\mathcal{G}_t]$. By Theorem 8.1 in Lipster and Shiryayev \cite{LipsterShiryayev77} and Theorem 5.7 (Kushner-FKK equation) in Xiong \cite{Xiong08}, we can get the state filtering equation for (\ref{afbsdde lq}):

\begin{equation}\label{xyhat lq}
\left\{
\begin{aligned}
d\hat{x}(t)=\ & [ A(t)\hat{x}(t)+\bar{A}(t)\hat{x}_{\delta}(t)-\sum_{i=1}^{2}B_i(t)R_i^{-1}(t)\mathcal{B}_i(t)]dt+[ C(t)\hat{x}(t)+\bar{C}(t)\hat{x}_{\delta}(t)\\
&-\sum_{i=1}^{2}D_i(t)R_i^{-1}(t)\mathcal{B}_i(t)]dW(t),\\
-d\hat{y}(t)=\ &\{E(t)\hat{x}(t)+F(t)\hat{y}(t)+G(t)\hat{z}(t)+\bar{G}(t)\hat{\bar{z}}(t)+\bar{F}(t)\mathbb{E}^{\mathcal{G}_t}[\hat{y}(t+\delta)]\\
&-\sum_{i=1}^{2}H_i(t)R_i^{-1}(t)\mathcal{B}_i(t)\}dt-\hat{z}(t)dW(t), \quad t\in[0,T],\\
\hat{x}(t)=\ &\xi(t),\quad t\in[-\delta,0],\\
\hat{y}(T)=\ &M_T\hat{x}(T),\ \hat{y}(t)=\hat{\varphi}(t),\quad t\in(T,T+\delta].\\
\end{aligned}
\right.
\end{equation}
where $\mathcal{B}_i(t)=B_i^{\tau}(t)\hat{q}_i(t)+D_i^{\tau}(t)\hat{k}_i(t)-H_i^\tau(t)\hat{p}_i(t)$. And the adjoint filtering equation for (\ref{ad lq}) satisfying

\begin{equation}\label{pqhat lq}
\left\{
\begin{aligned}
d\hat{p}_i(t)=\ &[F^\tau(t)\hat{p}_i(t)+\bar{F}^\tau(t-\delta)\hat{p}_i(t-\delta)-P_{i}(t)\hat{y}(t)]dt+[G^\tau(t)\hat{p}_i(t)-Q_{i}(t)\hat{z}(t)]dW(t)\\
-d\hat{q}_i(t)=\ &\{A^\tau(t)\hat{q}_i(t)+C^\tau(t)\hat{k}_i(t)-E^\tau(t)\hat{p}_i(t)+\mathbb{E}^{\mathcal{G}_t}[\bar{A}^\tau(t+\delta)\hat{q}_i(t+\delta)\\
&+\bar{C}^\tau(t+\delta)\hat{k}_i(t+\delta)]+O_i(t)\hat{x}(t)\}dt-\hat{k}_i(t)dW(t),\quad t\in[0,T],\\
\hat{p}_i(0)=&-N_iy(0),\ \hat{p}_i(t)=0,\quad t\in[-\delta,0),\\
\hat{q}_i(T)=&-M_T\hat{p}_i(T)+M_i\hat{x}(T),\ \hat{q}_i(t)=0,\quad t\in(T,T+\delta],\quad (i=1,2).
\end{aligned}
\right.
\end{equation}

From Theorem \ref{nec MP} and Theorem \ref{suf MP}, it is easy to know that $(u_1(\cdot),u_2(\cdot))$ is an equilibrium point for the above linear-quadratic game problem if and only if $(u_1(\cdot),u_2(\cdot))$ satisfies the expression of (\ref{ui lq}) with $(\hat{x},\hat{y},\hat{z},\hat{p}_i,\hat{q}_i,\hat{k}_i)(i=1,2)$ being the solution of the coupled triple dimensions filtering AFBSDDE (\ref{xyhat lq})-(\ref{pqhat lq}) (TFBSDDE). Then the existence and uniqueness of the equilibrium point is equivalent to the existence and uniqueness of the TFBSDDE. %Similar to HuangShi12, Yu12 and ChenWu11

However, the TFBSDDE (\ref{xyhat lq})-(\ref{pqhat lq}) is complicated, but, in some particular cases, we can use some transactions to relate it to a double dimensions filtering AFBSDDE, called DFBSDDE, such as the following result.

\textbf{H4.} The dimension of $x$ is equal to that of $y$: $n=m$, $\bar{G}(t)\equiv0$ and coefficients $B_i(t)=B_i,D_i(t)=D_i,H_i(t)=H_i$ are independent of time $t$ for any $i=1,2$.

\begin{theorem}
Under H4, we assume one of the following conditions holds true:

 (a) $D_1=D_2=H_1=H_2\equiv0$ and $B_iR^{-1}_iB^\tau_iS=SB_iR^{-1}_iB^\tau_i,\quad (i=1,2);$

 (b) $B_1=B_2=H_1=H_2\equiv0$ and $D_iR^{-1}_iD^\tau_iS=SD_iR^{-1}_iD^\tau_i,\quad (i=1,2)$;

 (c) $B_1=B_2=D_1=D_2\equiv0$ and $H_iR^{-1}_iH^\tau_iS=SH_iR^{-1}_iH^\tau_i,\quad (i=1,2)$,\\
 where $S^\tau=A(\cdot),\bar{A}(\cdot),C(\cdot),\bar{C}(\cdot)$,
 $E(\cdot),F(\cdot),\bar{F}(\cdot),G(\cdot),M_T,O_i(\cdot),P_i(\cdot),Q_i(\cdot),M_i,N_i$. Then $(u_1(\cdot)$,$u_2(\cdot))$ given by $(\ref{ui lq})$ is the unique Nash equilibrium point.\\
\end{theorem}

\textbf{Proof.}  We only proof (a). The same method can be used to get (b) and (c). From above discussion, we need to prove only that there exists a unique solution of the coupled TFBSDDE (\ref{xyhat lq})-(\ref{pqhat lq}). In the case that $D_1=D_2=H_1=H_2\equiv0$, it becomes

\begin{equation}\label{tafbsdde}
\left\{
\begin{aligned}
d\hat{x}(t)=\ & [ A(t)\hat{x}(t)+\bar{A}(t)\hat{x}_{\delta}(t)-\sum_{i=1}^{2}B_iR_i^{-1}B_i^\tau\hat{q}_i(t)]dt+[C(t)\hat{x}(t)+\bar{C}(t)\hat{x}_\delta(t)]dW(t),\\
-d\hat{y}(t)=\ &\{E(t)\hat{x}(t)+F(t)\hat{y}(t)+G(t)\hat{z}(t)+\bar{F}(t)\mathbb{E}^{\mathcal{G}_t}[\hat{y}(t+\delta)] \}dt-\hat{z}(t)dW(t),\\
d\hat{p}_i(t)=\ &[F^\tau(t)\hat{p}_i(t)+\bar{F}^\tau(t-\delta)\hat{p}_i(t-\delta)-P_{i}(t)\hat{y}(t)]dt+[G^\tau(t)\hat{p}_i(t)-Q_{i}(t)\hat{z}(t)]dW(t),\\
-d\hat{q}_i(t)=\ &\{A^\tau(t)\hat{q}_i(t)+C^\tau(t)\hat{k}_i(t)-E^\tau(t)\hat{p}_i(t)+\mathbb{E}^{\mathcal{G}_t}[\bar{A}^\tau(t+\delta)\hat{q}_i(t+\delta)\\
&+\bar{C}^\tau(t+\delta)\hat{k}_i(t+\delta)]+O_i(t)\hat{x}(t)\}dt-\hat{k}_i(t)dW(t),\quad t\in[0,T]\\
\hat{x}(t)=\ &\xi(t),\quad t\in[-\delta,0];\ \hat{y}(T)= M_T\hat{x}(T),\ \hat{y}(t)=\hat{\varphi}(t),\quad t\in(T,T+\delta],\\
\hat{p}_i(0)=&-N_iy(0),\ \hat{p}_i(t)=0,\quad t\in[-\delta,0),\\
\hat{q}_i(T)=&-M_T\hat{p}_i(T)+M_i\hat{x}(T),\ \hat{q}_i(t)=\hat{k}_i(t)=0,\quad t\in(T,T+\delta],\quad (i=1,2).
\end{aligned}
\right.
\end{equation}

Now we consider another DFBSDDE:

\begin{equation}\label{dafbsdde}
\left\{
\begin{aligned}
d\tilde{x}(t)=\ & [ A(t)\tilde{x}(t)+\bar{A}(t)\tilde{x}_{\delta}(t)-\tilde{q}(t)]dt+[C(t)\tilde{x}(t)+\bar{C}(t)\tilde{x}_\delta(t)]dW(t),\\
-d\tilde{y}(t)=\ &\{E(t)\tilde{x}(t)+F(t)\tilde{y}(t)+G(t)\tilde{z}(t)+\bar{F}(t)\mathbb{E}^{\mathcal{G}_t}[\tilde{y}(t+\delta)] \}dt-\tilde{z}(t)dW(t),\\
d\tilde{p}(t)=\ &[F^\tau(t)\tilde{p}(t)+\bar{F}^\tau(t-\delta)\tilde{p}(t-\delta)-\sum_{i=1}^2B_iR^{-1}_iB^\tau_iP_{i}(t)\tilde{y}(t)]dt+[G^\tau(t)\tilde{p}(t)\\
&-\sum_{i=1}^2B_iR^{-1}_iB^\tau_iQ_{i}(t)\tilde{z}(t)]dW(t),\\
-d\tilde{q}(t)=\ &\{A^\tau(t)\tilde{q}(t)+C^\tau(t)\tilde{k}(t)-E^\tau(t)\tilde{p}(t)+\bar{A}^{\tau}(t+\delta)\mathbb{E}^{\mathcal{G}_t}[\tilde{q}(t+\delta)]\\
&+\bar{C}^\tau(t+\delta)\mathbb{E}^{\mathcal{G}_t}[\tilde{k}(t+\delta)]+\sum_{i=1}^2B_iR^{-1}_iB^\tau_iO_{i}(t)\tilde{x}(t)\}dt-\tilde{k}(t)dW(t),\quad t\in[0,T]\\
\tilde{x}(t)=\ &\xi(t),\quad t\in[-\delta,0];\ \tilde{y}(T)=M_T\tilde{x}(T),\ \tilde{y}(t)=\hat{\varphi}(t),\quad t\in(T,T+\delta],\\
\tilde{p}(0)=&-\sum_{i=1}^2B_iR^{-1}_iB^\tau_iN_i\tilde{y}(0),\ \tilde{p}(t)=0,\quad t\in[-\delta,0),\\
\tilde{q}(T)=&\sum_{i=1}^2B_iR^{-1}_iB^\tau_iM_i\tilde{x}(T)-M_T\tilde{p}(T),\ \tilde{q}(t)=\tilde{k}(t)=0,\quad t\in(T,T+\delta].
\end{aligned}
\right.
\end{equation}

From the commutation relation between matrices, we notice that, if $(\hat{x},\hat{y},\hat{z},\hat{p}_i,\hat{q}_i,\hat{k}_i)(i=1,2)$ is a solution of (\ref{tafbsdde}), then $(\tilde{x},\tilde{y},\tilde{z},\hat{p},\tilde{q},\tilde{k})$ solves (\ref{dafbsdde}), where

\begin{equation*}
\left\{
\begin{aligned}
\tilde{x}(t)=&\hat{x}(t),\tilde{y}(t)=\hat{y}(t),\tilde{z}(t)=\hat{z}(t),\\
\tilde{p}(t)=&B_1R^{-1}_1B^\tau_1\hat{p}_1(t)+B_2R^{-1}_2B^\tau_2\hat{p}_2(t),\\
\tilde{q}(t)=&B_1R^{-1}_1B^\tau_1\hat{q}_1(t)+B_2R^{-1}_2B^\tau_2\hat{q}_2(t),\\
\tilde{k}(t)=&B_1R^{-1}_1B^\tau_1\hat{k}_1(t)+B_2R^{-1}_2B^\tau_2\hat{k}_2(t).\\
\end{aligned}
\right.
\end{equation*}

On the other hand, if $(\tilde{x},\tilde{y},\tilde{z},\hat{p},\tilde{q},\tilde{k})$ is a solution of  (\ref{dafbsdde}), we can let $\hat{x}(t)=\tilde{x}(t),\hat{y}(t)=\tilde{y}(t),\hat{z}(t)=\tilde{z}(t)$. Then we get $(\hat{p}_i(t),\hat{q}_i(t),\hat{k}_i(t))$ from the following filtering AFBSDDE:

\begin{equation*}
\left\{
\begin{aligned}
d\hat{p}_i(t)=\ &[F^\tau(t)\hat{p}_i(t)+\bar{F}^\tau(t-\delta)\tilde{p}_i(t-\delta)-P_{i}(t)\hat{y}(t)]dt+[G^\tau(t)\hat{p}_i(t)-Q_{i}(t)\hat{z}(t)]dW(t),\\
-d\hat{q}_i(t)=\ &\{A^\tau(t)\hat{q}_i(t)+C^\tau(t)\hat{k}_i(t)-E^\tau(t)\hat{p}_i(t)+\bar{A}^{\tau}(t+\delta)\mathbb{E}^{\mathcal{F}_t}[\hat{q}_i(t+\delta)]\\
&+\bar{C}^\tau(t+\delta)\mathbb{E}^{\mathcal{F}_t}[\hat{k}_i(t+\delta)]+\sum_{i=1}^2B_iR^{-1}_iB^\tau_iO_{i}(t)\hat{x}(t)\}dt-\hat{k}_i(t)dW(t),\quad t\in[0,T]\\
\hat{p}_i(0)=&-N_i\hat{y}(0),\ \hat{p}_i(t)=0,\quad t\in[-\delta,0),\\
\hat{q}_i(T)=&-M_T\hat{p}_i(T)+M_i\hat{x}(T),\ \hat{q}_i(t)=\hat{k}_i(t)=0,\quad t\in(T,T+\delta],\quad (i=1,2).
\end{aligned}
\right.
\end{equation*}

We let

\begin{equation}\label{pqk bar}
\left\{
\begin{aligned}
\bar{p}(t)=&B_1R^{-1}_1B^\tau_1\hat{p}_1(t)+B_2R^{-1}_2B^\tau_2\hat{p}_2(t),\\
\bar{q}(t)=&B_1R^{-1}_1B^\tau_1\hat{q}_1(t)+B_2R^{-1}_2B^\tau_2\hat{q}_2(t),\\
\bar{k}(t)=&B_1R^{-1}_1B^\tau_1\hat{k}_1(t)+B_2R^{-1}_2B^\tau_2\hat{k}_2(t),\\
\end{aligned}
\right.
\end{equation}

By It\^o's formula and the uniqueness result of the solution of the SDDE and ABSDE (see \cite{ChenWu10,PengYang09}) for fixed $(\hat{x}(\cdot),\hat{y}(\cdot),\hat{z}(\cdot))$,  we have

\begin{equation}\label{}
\left\{
\begin{aligned}
\tilde{p}(t)=\bar{p}(t)=&B_1R^{-1}_1B^\tau_1\hat{p}_1(t)+B_2R^{-1}_2B^\tau_2\hat{p}_2(t),\\
\tilde{q}(t)=\bar{q}(t)=&B_1R^{-1}_1B^\tau_1\hat{q}_1(t)+B_2R^{-1}_2B^\tau_2\hat{q}_2(t),\\
\tilde{k}(t)=\bar{k}(t)=&B_1R^{-1}_1B^\tau_1\hat{k}_1(t)+B_2R^{-1}_2B^\tau_2\hat{k}_2(t).\\
\end{aligned}
\right.
\end{equation}

This implies that the existence and uniqueness of (\ref{dafbsdde}) is equivalence to the existence and uniqueness of (\ref{tafbsdde}). According to the monotonic condition in \cite{ChenWu11,HuangShi12}, it is easy to check the DFBSDDE (\ref{dafbsdde}) satisfies the condition and it has a unique solution. So the TFBSDDE (\ref{tafbsdde}) admits a unique solution. We complete the proof.

\section{An example in finance}
This section is devoted to study a pension fund management problem under partial information with time-delayed surplus arising from the financial market, which naturally motivates the above theoretical research. The financial market is the Black-Scholes market, while the pension fund management framework comes from Federico \cite{Federico11}. To get close to reality, we  study this problem in the case when the performance criterion $J_i(v_1(\cdot),v_2(\cdot))$ is measured by a criterion involving \emph{risk}. If we interpret \emph{risk} in the sense of a \emph{convex risk measure}, it can be performed as a non-linear expectation called g-expectation, which can also be used to represent a non-linear human preference in behavioral economics. See \cite{FollmerSchied02,FrittelliGianin02,Peng97,Peng04} and recent articles \cite{HuiXiao12,AnOksendal12}. Now we introduce it in detail.

In the following, we only consider the 1-dimension case just for simplicity of notations. First, we give the definition of convex risk measure and its connection with g-expectation.

\begin{definition}[\cite{FollmerSchied02,FollmerSchield02a,FrittelliGianin02} ]\label{convex}
Let $\mathbb{F}$ be the family of all lower bounded $\mathcal{F}_T$-measurable random variables. A convex risk measure on $\mathbb{F}$ is a functional $\rho: \mathbb{F}\rightarrow \mathbb{R}$ such that

\noindent(a)(convexity) $\rho(\lambda X_1+(1-\lambda)X_2)\leq \lambda\rho(X_1)+(1-\lambda)\rho(X_2),\quad X_1,X_2\in \mathbb{F},\ \lambda\in(0,1),$

\noindent(b)(monotonicity) if $X_1\leq X_2$ a.e., then $\rho(X_1)\geq \rho(X_2),\quad X_1,X_2\in\mathbb{F},$

\noindent(c)(translation invariance) $\rho(X+m)=\rho(X)-m,\quad X\in\mathbb{F},m\in\mathbb{R}.$
\end{definition}

The convex risk measure is a useful tool widely applied in the measurement of financial positions. The property (a) in Definition \ref{convex} represents a non-linearity that illustrates the better choice for the diversified investments; (b) means that if portfolio $X_2$ is better than $X_1$ under almost all scenarios, then the risk of $X_2$ should be less than the risk of $X_1$; (c) implies that the addition of a sure amount of capital reduces the risk by the same amount. It is also a generalization of the concept of \emph{coherent risk measure} in \cite{ArtznerDelbaen99}.

Consider the following BSDE:

\begin{equation}\label{y example}
\left\{
\begin{aligned}
-dy(t)&=g(t,z(t))dt-z(t)dW(t),\\
y(T)&=\xi.\\
\end{aligned}
\right.
\end{equation}

Under certain assumptions, (\ref{y example}) exists a unique solution $(y(\cdot),z(\cdot))$. If we also set $g(\cdot,0)\equiv0$, we can make the definition as follows.

\begin{definition}[\cite{Peng97,Peng04} ]\label{gexpectation}
For each $\xi\in\mathcal{F}_T$, we call

\begin{equation*}
\mathcal{E}_g(\xi)\triangleq y(0)
\end{equation*}
the generalized expectation (g-expectation) of $\xi$ related to $g$.
\end{definition}

We can know that the expectation $\mathbb{E}$ is a linear expectation and does not express peoples's \emph{preferences} or criterion involving \emph{risk}, and the map $\xi\rightarrow \mathcal{E}_g(\xi)$ includes all the properties that $\mathbb{E}$ has, except the linearity. It is obvious that when $g(\cdot)=0$, $\mathcal{E}_g$ is reduced to the classical expectation $\mathbb{E}$.

Now, we give a connection between the convex risk measure and the g-expectation.
\begin{definition}\label{connection}
The risk $\rho(\xi)$ of the random variable $\xi\in\mathcal{L}^2_\mathcal{F}(\Omega;\mathbb{R})$ ($\xi$ can be regarded as a financial position in the financial market) is defined by

\begin{equation*}
\rho(\xi)\triangleq\mathcal{E}_g[-\xi]=y(0),
\end{equation*}
where $\mathcal{E}_g[\cdot]$ is defined in the Definition \ref{gexpectation}, but with $\xi$ replaced by $-\xi$.
\end{definition}

Assuming that there are two asset in the financial market for the pension fund managers to invest:

\begin{equation*}
\left\{
\begin{aligned}
dS_0(t)&= r(t)S_0(t)dt,\\
dS_1(t)&= \mu(t)S_1(t)dt+\sigma(t)S_1(t)dW(t),\\
S_0(0)&= 1,\  S_1(0)> 0,
\end{aligned}
\right.
\end{equation*}
where $S_1(\cdot)$ is a risky finance asset price and $S_0(\cdot)$ is one risk-free asset price. $\mu(\cdot)$ is an appreciation rate of the asset process, and the $\sigma(\cdot)$ is the volatility coefficients. We assume that $\mu(\cdot)$, $r(\cdot)$ and $\sigma(\cdot)$ are deterministic bounded coefficients, and $\sigma^{-1}(\cdot)$ is bounded.

Suppose that there are two pension fund managers (players) working together to invest the risk-free and risky assets. In real financial market, it is reasonable for the investors to make decisions based on the historical price of the risky asset $S_1(\cdot)$. So the observable filtration can be set as $\mathcal{G}_t=\sigma\{S_1(s)|0\leq s\leq t\}$, and it is clear that $\mathcal{G}_t=\mathcal{F}_t^{W}=\sigma\{W(s)|0\leq s\leq t\}$. The pension fund wealth $x(\cdot)$ can be modeled by

\begin{equation}\label{x example}
\left\{
\begin{aligned}
dx(t)=\ &(r(t)x(t)+(\mu(t)-r(t))\pi(t)-\alpha(x(t)-x(t-\delta))-c_1(t)-c_2(t))dt+\pi(t)\sigma(t)dW(t)\\
&+\bar{\sigma}(t)d\bar{W}(t), \\
x(0)=\ &x_0>0,\ x(t)=0,\quad t\in[-\delta,0).\\
\end{aligned}
\right.
\end{equation}
Here we denote by $\pi(t)$ the amount of portfolio invested in the risky asset at time $t$, and $\alpha(x(t)-x(t-\delta))$ represents the surplus premium to fund members or their capital transfusions depending on the performance of fund growth during the past period with parameter $\alpha>0$ (see e.g. \cite{WuWang15,ShiWang15}). Meanwhile, there is an instantaneous consumption rate $c_i(t)$ for manager $i(i=1,2)$. We assume that the value of $x(\cdot)$ is not only affected by the risky asset, but also by some practical phenomena like the physical inaccessibility of some economic parameters, inaccuracies in measurement, discreteness of account information, insider trading or the information asymmetry between the manager and investors, etc (see e.g. \cite{Lakner95,HuangWangWu10}). Thus we set $\bar{\sigma}(\cdot)$ be the instantaneous volatility, $\mathcal{F}_t^{\bar{W}}$ represents the unobservable filtration generated by $\bar{W}(\cdot)$, $x(t)$ be adapted to the filtration $\mathcal{F}_t$ generated by Brownian motion $(W(\cdot),\bar{W}(\cdot))$, and the control processes $c_i(t)\ (i=1,2)$ be adapted to the observation filtration $\mathcal{G}_t\subseteq\mathcal{F}_t$.

The controlled processes $c_i(\cdot)\ (i=1,2)$ is called \emph{admissible} for manager $i$ if $c_i(t)>0$ is adapted to the filtration $\mathcal{G}_t$ at time $t$, $c_i(t)\in\mathbb{L}^2(0,T;\mathbb{R})$, and the family of admissible control $(c_1(\cdot),c_2(\cdot))$ is denoted by $\mathcal{C}_1\times\mathcal{C}_2$.

We assume that the insurance company hopes more terminal capital under less risk and more consumption $c_i(\cdot)$. According to the Definition \ref{convex} and \ref{connection}, we can define the cost functional as
\begin{equation}\label{Jg}
J^g_i(c_1(\cdot),c_2(\cdot))=-\mathcal{E}_g[-x(T)]+\mathbb{E}\int_0^Te^{-\beta t}L_i\frac{c_i(t)^\gamma}{\gamma}dt,\qquad i=1,2
\end{equation}
where $L_i$ is a positive constant, $\beta$ is a discount factor and $1-\gamma\in(0,1)$ is a constant called the Arrow-Pratt index of risk aversion.

Then our problem is naturally to find an equilibrium point $(c_1^*(\cdot),c_2^*(\cdot))\in\mathcal{C}_1\times\mathcal{C}_2$ such that

\begin{equation*}
\left\{
\begin{aligned}
J_1^g(c_1^*(\cdot),c_2^*(\cdot))&=\sup\limits_{c_1\in\mathcal{C}_1}J_1^g(c_1(\cdot),c   _2^*(\cdot)),\\
J_2^g(c_1^*(\cdot),c_2^*(\cdot))&=\sup\limits_{c_2\in\mathcal{C}_2}J_2^g(c_1^*(\cdot),c_2(\cdot)).\\
\end{aligned}
\right.
\end{equation*}

In the following, we set $g(\cdot)$ be a linear form as $g(\cdot,z(\cdot))=g(\cdot)z(\cdot)$ where $g(\cdot)$ is deterministic bounded coefficient. Then our problem can be reformulated as

\begin{equation}\label{example xy}
\left\{
\begin{aligned}
dx(t)=\ &(r(t)x(t)+(\mu(t)-r(t))\pi(t)-\alpha(x(t)-x(t-\delta))-c_1(t)-c_2(t))dt+\pi(t)\sigma(t)dW(t)\\
&+\bar{\sigma}(t)d\bar{W}(t),\\
-dy(t)=\ &g(t)z(t)dt-z(t)dW(t),\quad t\in[0,T],\\
x(0)=\ & x_0,\ \ x(t)=0,\quad t\in[-\delta,0),\\
y(T)=\ &-x(T),
\end{aligned}
\right.
\end{equation}
and

\begin{equation}\label{example Jx}
J^g_i(c_1(\cdot),c_2(\cdot))=\mathbb{E}\int_0^Te^{-\beta t}L_i\frac{c_i(t)^\gamma}{\gamma}dt-y(0),\qquad i=1,2.
\end{equation}

Now we will apply the theoretical results obtained in Section \ref{MP} to solve the above game problem. The Hamiltonian function is in the form of

$$
\begin{aligned}
H_i&(t,x(t),y(t),z(t),x_\delta(t),c_1(t),c_2(t);p(t),q_i(t),k_i(t),\bar{k}_i(t))=q_i(t)[r(t)x(t)+(\mu(t)-r(t))\pi(t)\\
&-\alpha(x(t)-x(t-\delta))-c_1(t)-c_2(t)]+k_i(t)\pi(t)\sigma(t)+\bar{k}_i(t)\bar{\sigma}(t)-p(t)g(t)z(t)+e^{-\beta t}L_i\frac{c_i(t)^\gamma}{\gamma},
\end{aligned}
$$

where the adjoint process satisfies

\begin{equation*}
\left\{
\begin{aligned}
dp(t)=\ &g(t)p(t)dW(t),\\
-dq_i(t)=\ &\{(r(t)-\alpha)q_i(t)+\alpha\mathbb{E}^{\mathcal{F}_t}[q_i(t+\delta)]\}dt-k_i(t)dW(t)-\bar{k}_i(t)d\bar{W}(t),\quad t\in[0,T],\\
p(0)=\ &1,\\
q_i(T)=\ &p(T),\ q_i(t)=k_i(t)=\bar{k}_i(t)=0,\quad t\in(T,T+\delta],\quad (i=1,2).
\end{aligned}
\right.
\end{equation*}

Then we use the necessary maximum principle (Theorem \ref{nec MP}) to find a candidate equilibrium point:

\begin{equation}\label{c1c2}
\begin{aligned}
  c_1^*(t)=(L_1^{-1}e^{rt}\hat{q}_1(t))^{\frac{1}{\gamma-1}},\\
  c_2^*(t)=(L_2^{-1}e^{rt}\hat{q}_2(t))^{\frac{1}{\gamma-1}}.
  \end{aligned}
\end{equation}
where $\hat{q}_i(t)=\mathbb{E}[q_i(t)|\mathcal{G}_t]\ (i=1,2)$.

Now we have to deal with $\hat{q}_i(t)$, the optimal filtering of $q_i(t)$ on the observation $\mathcal{G}_t$. We also set $\hat{p}_i(t)=\mathbb{E}[p_i(t)|\mathcal{G}_t]$. Note that

$$
\begin{aligned}
\mathbb{E}\{\mathbb{E}[q_i(t+\delta)|\mathcal{F}_t]|\mathcal{G}_t\}=\mathbb{E}[q_i(t+\delta)|\mathcal{G}_t]=\mathbb{E}\{\mathbb{E}[q_i(t+\delta)|\mathcal{G}_{t+\delta}]|\mathcal{G}_t\}
=\mathbb{E}[\hat{q}_i(t+\delta)|\mathcal{G}_t].
\end{aligned}
$$

Then by Theorem 8.1 in \cite{Xiong08}, we have

\begin{equation}\label{pqhat example}
\left\{
\begin{aligned}
d\hat{p}(t)=\ &g(t)\hat{p}(t)dW(t),\\
-d\hat{q}_i(t)=\ &\{(r(t)-\alpha)\hat{q}_i(t)+\alpha\mathbb{E}^{\mathcal{G}_t}[\hat{q}_i(t+\delta)]\}dt-\hat{k}_i(t)dW(t),\quad t\in[0,T],\\
\hat{p}(0)=\ &1,\\
\hat{q}_i(T)=\ &\hat{p}(T),\ \hat{q}_i(t)=\hat{k}_i(t)=0,\quad t\in(T,T+\delta],\quad (i=1,2).
\end{aligned}
\right.
\end{equation}

From (\ref{pqhat example}), we can derive the explicit expression of $\hat{p}_i(t)$ as

\begin{equation*}
  \hat{p}(t)=\hat{p}_i(t)=\exp\{\int_0^tg(s)dW(s)-\frac{1}{2}\int_0^tg^2(s)ds\}>0,\quad t\in[0,T],
\end{equation*}
which is an $\mathcal{G}_t$-exponential martingale.

By Theorem 5.1 in \cite{PengYang09}, we can prove $\hat{q}_i(t)\geq 0,\ t\in[0,T]$. Thus $c^*_i(t)>0$ for all $t\in[0,T]$. Next, we will solve the anticipated BSDE of $\hat{q}_i(t)$ recursively. This method can also be found in \cite{Yu12,MenoukeuPamen11}.

(1) When $t\in[T-\delta,T]$, the ABSDE in (\ref{pqhat example}) becomes a standard BSDE (without anticipation):

$$
\hat{q}_i(t)=\hat{p}(T)+\int_t^T (r(s)-\alpha)\hat{q}_i(s)ds-\int_t^T\hat{k}_i(s)dW(s),\quad t\in[T-\delta,T].
$$

Obviously, we have
$$
\hat{q}_i(t)=\exp\{\int_t^T(r(s)-\alpha)ds\}\mathbb{E}^{\mathcal{G}_t}[\hat{p}(T)]=\exp\{\int_t^T(r(s)-\alpha)ds\}\hat{p}(t),\quad t\in[T-\delta,T].
$$

From Proposition 5.3 in \cite{KarouiPengQuenez97}, $(\hat{q}_i(t),\hat{k}_i(t))$ is Malliavin differentiable and $\{D_t\hat{q}_i(t);T-\delta\leq t\leq T\}$ provides a version of  $\{\hat{k}_i(t);T-\delta\leq t\leq T\}$, i.e.

$$
\hat{k}_i(t)=D_t\hat{q}_i(t)=\exp\{\int_t^T(r(s)-\alpha)ds\}D_t\hat{p}(t),\quad t\in[T-\delta,T].
$$

(2) If we have solved ABSDE (\ref{pqhat example}) on the interval $[T-n\delta,T-(n-1)\delta](n=1,2,...)$, and the solution $\{(\hat{q}_i(t),\hat{k}_i(t));T-n\delta\leq t\leq T-(n-1)\delta\}$ is Malliavin differentiable, then we continue to consider the solvability on the next interval $[T-(n+1)\delta,T-n\delta]$, where we can rewrite the ABSDE (\ref{pqhat example}) as follows:

$$
\hat{q}_i(t)=\hat{q}_i(T-n\delta)+\int_t^{T-n\delta}\{(r(s)-\alpha)\hat{q}_i(s)+\alpha\mathbb{E}^{\mathcal{G}_s}[\hat{q}_i(s+\delta)]\}ds-\int_t^{T-n\delta}\hat{k}_i(s)dW(s).
$$

We note that $\{(\hat{q}_i(s+\delta),\hat{k}_i(s+\delta));t\leq s\leq T-n\delta\}$ has been solved and is Malliavin differentiable. So the same discussion lead to $\{(\hat{q}_i(t),\hat{k}_i(t));T-(n+1)\delta\leq t\leq T-n\delta\}$ is Malliavin differentiable, and

$$
\hat{q}_i(t)=\exp\{\int_t^{T-n\delta}(r(s)-\alpha)ds\}\mathbb{E}^{\mathcal{G}_t}[\hat{q}_i(T-n\delta)]+\alpha\int_t^{T-n\delta}\exp\{\int_t^s(r(\eta)-\alpha)d\eta\}\mathbb{E}^{\mathcal{G}_t}[\hat{q}_i(s+\delta)]ds,
$$

$$
\hat{k}_i(t)=\exp\{\int_t^{T-n\delta}(r(s)-\alpha)ds\}\mathbb{E}^{\mathcal{G}_t}[D_t\hat{q}_i(T-n\delta)]+\alpha\int_t^{T-n\delta}\exp\{\int_t^s(r(\eta)-\alpha)d\eta\}\mathbb{E}^{\mathcal{G}_t}[D_t\hat{q}_i(s+\delta)]ds,
$$
for any $t\in[T-(n+1)\delta,T-n\delta]$, $i=1,2$.

We notice that all the condition in the verification theorem (Theorem \ref{suf MP}) are satisfied, then Theorem \ref{suf MP} implies that $(c_1^*(\cdot),c_2^*(\cdot))$ given by (\ref{c1c2}) is an equilibrium point.

\begin{pro}
  The investment problem (\ref{x example})-(\ref{Jg}) admits an equilibrium point $(c_1^*(\cdot),c_2^*(\cdot))$ which is defined by (\ref{c1c2}).
\end{pro}

\end{document}